\documentclass[11pt,showlabels]{article1}
\usepackage{amssymb}
\usepackage{epsfig}
\newtheorem{theorem}{Theorem}%[section]
\newtheorem{assumption}{Assumption}%[section]
\newtheorem{lemma}{Lemma}%[section]
\newtheorem{remark}{Remark}%[section]
\newtheorem{definition}{Definition}%[section]
\newtheorem{corollary}{Corollary}%[section]

\title{\bf Linear Quadratic Optimal Control and Stabilization for Discrete-time Markov Jump Linear Systems
 \thanks{This work is supported by the National Natural Science
Foundation of China (Nos. 61473134, 61573220, 61120106011, 61573221). $^{*}$Corresponding author: Huanshui Zhang. Email: hszhang@sdu.edu.cn}
}

\author{Chunyan Han$^a$,\ Hongdan Li$^b$, \ Wei Wang$^{b}$, \ Huanshui Zhang$^{b,*}$
\ \\
\\
\ \ \ $^a$ School of Electrical Engineering, University of Jinan, \\Jinan Shandong 250022, China
\\
$^b$ School of Control Science and Engineering, Shandong
University, \\Jinan Shandong 250061, China}
\begin{document}
\baselineskip 16pt
\date{}
%\begin{titlepage}
  \maketitle
\begin{abstract}
%This paper investigates the finite and infinite-horizon optimal control problems for the Markov jump linear systems. The general case for the finite-horizon optimal controller is considered, where the input weighting matrix in the performance index is just required to be positive semidefinite. The necessary and sufficient condition for the existence of the optimal finite-horizon controller is given, and the optimal controller is derived explicitly from a set of coupled difference Riccati equations (CDRE). For the case of infinite-horizon case, the primary concerning problem is the mean square stabilizability and optimality of optimal controller. A necessary and sufficient condition is established for the existence of the mean square stabilizable solution to the optimal controller for the first time, and the stabilization condition is guaranteed by the positiveness of the associated coupled algebraic Riccati equation (CARE) which is easy to be shown. Meanwhile, the optimality of the solution is guaranteed. The optimal stabilizing solution to the infinite-horizon case is given by deriving the corresponding CARE. By a product, the sufficient conditions for the convergence of the CDRE to the CARE is supplied.

This paper mainly investigates the optimal control and stabilization problems for linear discrete-time Markov jump systems.  The general case for the finite-horizon optimal controller is considered, where the input weighting matrix in the performance index is just required to be positive semi-definite.  The necessary and sufficient condition for the existence of the optimal controller in finite-horizon is given explicitly from a set of coupled difference Riccati equations (CDRE). One of the key techniques is to solve the forward and backward stochastic difference equation (FDSDE) which is obtained by the maximum principle. As to the infinite-horizon case, we  establish the necessary and sufficient condition to  stabilize the Markov jump linear system in the mean square sense. It is shown that the Markov jump linear system is stabilizable under the optimal controller if and only if the associated couple algebraic Riccati equation (CARE) has a unique positive solution. Meanwhile, the optimal controller and optimal cost function in infinite-horizon case are expressed explicitly.
\bigskip

\end{abstract}

\pagestyle{plain} \setcounter{page}{1}
\section{Introduction}
Over the last decades, there has been a steadily rising level of activity with linear systems subject to abrupt changes in their structures. The case in which the jumps are modeled by a Markov chain is referred to as Markov jump linear systems (MJLS) and has been receiving lately a great deal of attention in the literature. Applications of these models can be found, for instance, in robotics tracking and estimation, communication networks, flight systems, etc. We can mention the books and the references therein for a general overview on the control and filter problems for MJLS \cite{CFM05}, \cite{Mar1990}.

The study of this discrete-time Markovian jump linear quadratic (JLQ) control problem can be traced back (at least) to the work of Blair and Sworder \cite{BS75} for the finite-time horizon case. Sworder used the dynamic programming method to obtain his result. Birdwell et al. \cite{BCA86} examined the case where matrix A is not dependent on the form process. For the infinite time horizon version of this control problem, necessary and sufficient conditions for existence of the optimal steady-state JLQ controller were provided in \cite{CWC86}, and sufficient conditions for these steady-state control laws to stabilize the controlled system were given. These conditions were not easily tested, however, since they required the simultaneous solution of coupled matrix equations containing infinite sums. Later in \cite{JC88}, new and refined definitions of the controllability and observability of discrete-time MJLS were developed. Algebraic test for these concepts were also given, and the existence of optimal steady-state JLQ controllers was guaranteed by the absolute controllability. Under the precondition of the existence of a steady-state controller and the finiteness of the optimal expected cost, the absolute observability guaranteed the stability of the controlled system. This result avoided the awkward need of finding constant (non-optimal) control laws in order to check for steady-state convergence, as in \cite{CWC86}. However, this result involved only sufficient conditions and finiteness of the optimal expected cost needed to be checked. In \cite{JC89} and \cite{JC1989}, the JLQ problem for systems with $x$- and $u$-dependent form transition probabilities were considered. In \cite{AFJ95}, a necessary and sufficient condition was presented for the existence of a positive-semidefinite solution of the coupled algebraic Riccati-like equation occurring in the infinite horizon JLQ problems. However, the existence of the optimal controller and the stability of the closed-loop system were not discussed in the paper. In \cite{CFM05}, both finite horizon and infinite horizon optimal controls were considered, and the optimal controllers derived from a set of coupled difference Riccati equations(CDRE) for the former problem, and the mean square stabilizing solution of a set of coupled algebraic Riccati equations (CARE) for the latter problem. Sufficient existence condition for the finite-horizon optimal controller was guaranteed by the positive definite of the input weighting matrix, and the mean square stabilizing solution to the infinite-horizon optimal control problem was derived under the assumption that the stabilizing solution to CARE was existed.
Subsequently, the results developed in \cite{CFM05} were extended to solve the constrained quadratic control problems \cite{CABM99} and the finite-horizon output feedback quadratic optimal control problems \cite{CT03}, respectively.
%In [8], the constrained quadratic control problem of discrete-time MJLS has been considered. A state feedback mean-square stabilizing controller was obtained by solving a convex optimization problem via linear matrix inequality method. In [9], the finite-horizon output feedback quadratic optimal control problem of discrete-time MJLS was considered. A separation for this problem was established, and an optimal controller was obtained from two coupled Riccati difference equations under the assumption of positive definite of the input weighting matrix appearing in the performance index.
In \cite{CV01}, a new detectability concept (weak detectability) for discrete-time MJLS was presented, and the new concept supplied a sufficient condition for the mean square stable for the infinite-horizon linear quadratic controlled system. The concept of weak detectable retrieved the idea that each non-observed state corresponded to stable modes of the system. It has been shown in \cite{CV01} that mean square detectability ensured weak detectability. Latter, Costa and do Val \cite{CV02} summarized the available results and gave a proposition on the mean square stabilizable of the system. Under the assumption that the system was weak detectable, the system was mean square stabilizable if and only if there existed a positive semi-definite solution to the CARE. And a method for seeking the stabilizing solution of the CARE was supplied in \cite{CV02}. Further in \cite{CF15}, it studied the state-feedback JLQ control problem for discrete-time Markov jump linear systems considering the case in which the Markov chain takes values in a general Borel space. It was shown that the solution of the JLQ optimal control problem was obtained in terms of the positive semi-definite solution of M-coupled algebraic Riccati equations. It was obtained sufficient conditions, based on the concept of stochastic stabilizability and stochastic detectability, for the existence and uniqueness of this positive semi-definite solution \cite{CFM05}. Meanwhile, the output feedback JLQ optimal control problem for this type of systems was studied in \cite{CF16}.

Continuous-time version of the jump linear quadratic control problem were solved for finite-time horizons in Sworder \cite{Swor1969} and Wonham \cite{Wonh1971}. Sworder used a stochastic maximum principle to obtain his result, and Wonham used dynamic programming. Wonham also solved the infinite time horizon version of this control problem, and derived a set of sufficient conditions for the existence of a unique, finite steady-state solution. Mariton \cite{MB85} considered a discount cost version of the problem, where a controller that ensures stability in all forms was obtained. However, the paper used the unstated hypothesis that the diagonal entries of the generator of the jump process are equal. A discussion about this result appears in \cite{Hopk1986}. Mariton and Bertrand \cite{MB1985} also considered an output feedback version of the JLQ problem. Necessary optimality conditions are derived and two computational algorithms proposed. However, it was not possible to demonstrate their convergence. In \cite{JC90}, a necessary and sufficient condition for stochastic stabilizability of the system was provided, where the concept of stochastic stabilizability means the boundedness of the infinite-horizon performance index. Under the prerequisite that the optimal solution of the infinite-horizon JLQ problem existed, a necessary and sufficient stabilizing condition for the JLQ solution was provided in \cite{VC05} via the definition of weak detectable. In \cite{TW12}, a stochastic maximum principle for the finite-horizon optimal control problems of the continuous-time forward-backward Markovian regime-switching system was provided. The control system was described by forward-backward stochastic differential equations and modulated by continuous-time, finite-state Markov chains. The necessary and sufficient conditions for the optimal control was obtained.

In summary, the aforementioned works have supplied good results for the advances of the finite and infinite-horizon optimal control theories. As regards the discrete-time finite-horizon optimal control problems for the MJLS, the sufficient conditions for the existence of the optimal controller is guaranteed by the positiveness of the input penalty matrix $R$ or the positiveness of a set of matrix expressions, but no necessary conditions were presented. As for the discrete-time infinite-horizon optimal stabilization control problem, the necessary conditions \cite{CFM05}, \cite{CABM99},\cite{CT03}, sufficient conditions \cite{JC88}, \cite{CV01}, and necessary and sufficient conditions \cite{CWC86}, \cite{CV02} for the existence of the optimal stabilization controllers were provided. It need to point out that the necessary and sufficient conditions supplied in \cite{CWC86} were not easily tested, however, since they required the simultaneous solution of coupled matrix equations containing infinite sums. And the conditions developed in \cite{CV02} was based on the concept of weak detectability, it was not easy to test. To find a necessary and sufficient conditions, which is easy to check, for mean square stabilizable of the system in the linear optimal control frame is still an interesting problem. In the most recent works, Zhang et. al., \cite{ZWL12}-\cite{LZ16} considered the linear quadratic regulation (LQR) and stabilizaiton problem for the multiplicative noise systems with input delays. The necessary and sufficient condition for the existence of the finite-time LQR controller was established, and an explicit solution was given based on the solving forward and backward stochastic deferential/difference equations (FBSDEs)  which are from the maximum principle (MP).
%It was also shown that the system was stabilizable in the mean square sense if and only if the corresponding generalized algebraic Riccati equation(GARE) has a unique positive definite solution.
%The results based on a new technical tool (the stochastic maximum principle) and the introduction of a new Lyapunov function, in which a delayed forward-backward stochastic difference equation (D-FBSDE) needed to solve.
Inspired on the results developed in\cite{ZWL12}-\cite{LZ16}, we will propose a new approach to the LQR and stabilization problems for the MJLS based on MP. Compared to the multiplicative noise systems, the jumping parameter systems become more complicated since the correlation of the jumping parameters at adjoining time. We assume that the state variable and the jump parameters are available to the controller. In the second part of this paper we address the finite-horizon LQR for the discrete-time MJLS, where the input penalty matrix R is just required to be semi-definite positive. This relaxes the constraint imposed in the control problems greatly, which form the first innovation of this paper. We first extend the stochastic maximum principle \cite{ZWL12} to the jumping parameter systems, and develop a new forward-backward Markov jumping difference equation(FBMJDE). By solving the FBMJDE, a necessary and sufficient condition for the existence of the optimal controller is given, and an explicit analytical expression is given for the optimal controller. In the third part, a necessary and sufficient condition for the stabilization solution to the infinite-horizon optimal control problem is provided and the optimal constant gain controller is expressed explicitly, and the finite value of the infinite-horizon performance is given. Also, a special case is considered in Corollary 1. Compared with the existed results \cite{CFM05}, \cite{CV01}, \cite{CV02}, the stabilization condition in Corollary 1 is easy to test, since no precondition needs to test and it just requires to determine the existence of a positive definite solution to a set of CARE.
% By a product, the convergence of the CDRE concerning with the finite-horizon control are also shown, which converge to the set of CARE concerning the infinite-horizon control problem.
 The stochastic maximum principle and the explicit relationship between the optimal costate and the systems state explored in this paper play an important role in the derivation of the results.

Notations: Throughout this paper, ${R}^n$ denotes the $n$-dimensional Euclidean space,
$R^{m\times n}$ denotes the norm bounded linear space of all
$m\times n$ matrices.
For $L\in R^{n\times n}$, $L'$ stands for the transpose of $L$. As usual, $L\geq 0 (L>0)$
will mean that the symmetric matrix $L\in R^{n\times n}$ is positive semi-definite
(positive definite), respectively.

%%%%%%%%%%%%%%%%%%%%%%%%%%%%%%%%%%%%%%%%%%%%%%%%%%%%%%%%%%%%%%%%%%%%%%%%%%%%%%%%%%%%%%%%%%%%%%%%%%%%%%%%%%%%%%%%%%%%%%%%%%%%%%%%%%%%%%%%%%%%%%%%%%%%%%%%%%%%%%
\section{Finite-Horizon LQR for MJLS}

%%%%%%%%%%%%%%%%%%%%%%%%%%%%%%%%%%%%%%%%%%%%%%%%%%%%%%%%%%%%%%%%%%%%%%%%%%%%%%%%%%%%%%%%%%%%%%%%%%%%%%%%%%%%%%%%%%%%%%%%%%%%%%%%%%%%%%%%%%%%%%%%%%%%%%%%%%%%%%
\subsection{Problem Statement}

\setcounter{equation}{0}
We consider in this paper the finite horizon
optimal control problem for the Markov jump linear system (MJLS) when the state variable $x(k)$ and the jump variable $\theta(k)$ are available to the controller. On the stochastic basis $(\Omega,{\cal{G}}, {\cal{G}}_k,\mbox{P})$, consider the following MJLS
\begin{eqnarray}
x(k+1)&=&A_{\theta(k)}(k)x(k)+B_{\theta(k)}(k)u(k),\label{ff1}
\end{eqnarray}
where $x(k)\in {\mbox{R}}^n$ is the state, $u(k)\in {\mbox{R}}^m$  is the input control. $\theta(k)$ is a discrete-time Markov chain with finite state space $\{1,2,\cdots,L\}$ and transition probability $\lambda_{i,j}=\mbox{P}(\theta(k+1)=j|\theta(k)=i)(i,j=1,2,\cdots,L)$. We set $\pi_{i}(k)=\mbox{P}(\theta(k)=i)(i=1,2,\cdots,L)$, while $A_{i}(k), B_{i}(k)({i}=1,\cdots,L)$ are matrices of appropriate dimensions. The initial value $x_0$ is known. We assume that $\theta(k)$ is independent of $x_0$.

The quadratic cost
associated to system (\ref{ff1}) with admissible control law $u=(u(0),\cdots,u(N))$ is given by
\begin{eqnarray}
J_N&=&\mbox{E}[\sum_{k=0}^Nx(k)'Q_{\theta(k)}(k)x(k)+\sum_{k=0}^Nu(k)'R_{\theta(k)}(k)u(k)\nonumber\\
&&+x(N+1)'P_{\theta(N+1)}(N+1)x(N+1)],\label{f2}
\end{eqnarray}
where $N>0$ is an integer, $x(N+1)$ is the terminal state, $P_{j}(N+1)(j=1,\cdots,L)$ reflects the penalty on the terminal state, the matrix functions $R_{i}(k)\geq0(i=1,\cdots,L)$ and $Q_{i}(k)\geq 0(i=1,\cdots,L)$. The controller is required to obey the causality constraint, i.e., $u(k)$ must be in the form of
\begin{eqnarray}
u(k)=f_k(\theta(k),x(k),\cdots,x(0),u(k-1),\cdots,u(0))\nonumber
\end{eqnarray}
for some function $f_k(.)$. It means that $u(k)$ must be ${\cal{G}}_k$-measurable, where ${\cal{G}}_k=\{\theta(t);t=0,\cdots,k\}$. So the linear quadratic regulation (LQR) problem for Markov jumping parameter system can be stated as follows:

\emph{Problem 1}: Find a ${\cal{G}}_k$-measurable $u(k)$ such that (\ref{f2}) is minimized, subject to (\ref{ff1}).

\begin{remark}
 For brevity, we will omit the time steps in the systems matrices and the penalty matrices in the following discussions. That is denoting $A_{\theta(k)}(k), B_{\theta(k)}(k)$, $Q_{\theta(k)}(k)$, and $R_{\theta(k)}(k)$ as $A_{\theta(k)}, B_{\theta(k)}$, $Q_{\theta(k)}$, and $R_{\theta(k)}$, respectively. This will not affect the final results.
 \end{remark}
 \begin{remark}
 For finite-horizon optimal control of discrete-time MJLS, some results have been obtained. When the system is described by linear difference equations and when the penalty matrix $R_{\theta(k)}$ in the performance index is positive definite, the formalism of dynamic programming may be applied to advantage \cite{BS75}, \cite{CFM05}. However, only sufficient conditions to guarantee the existence of the optimal controller were given in \cite{BS75}, \cite{CFM05}. It is not possible to demonstrate a necessary and sufficient condition subject to the general case of $R_{\theta(k)}\geq0$. So in this paper, we consider the general case that $R_{\theta(k)}\geq0$. One object of this paper is to extend the work on linear system with white Gaussian noise to MJLS. To do this it will be expedient to derive an algorithm similar to the stochastic maximum principle developed in \cite{ZWL12}. The use of maximum principle in MJLS may supply a necessary and sufficient condition for the existence of the finite-horizon optimal controller in the general case.
 \end{remark}

%In next section, we will propose an extended version of the maximum principle. Different from the stochastic maximum principle for the system with multiplicative noises, the maximum principle developed in the following is suitable for the MJLS and provide a ${\cal{G}}_{k}$-measurable  controller $u(k)$. By employing the new proposed maximum principle, an explicit solution to the LQR problem will be obtained without state augmentation.

%%%%%%%%%%%%%%%%%%%%%%%%%%%%%%%%%%%%%%%%%%%%%%%%%%%%%%%%%%%%%%%%%%%%%%%%%%%%%%%%%%%%%%%%%%%%%%%%%%%%%%%%%%%%%%%%%%%%%%%%%%%%%%%%%%%%%%%%%%%%%%%%%%%%%%%%%%%%%%%
\subsection{Solution to the Finite-horizon LQR}

In the next, we will derive the optimal control by employing the stochastic maximum principle, where the necessary and sufficient condition for the existence of the optimal controller is proposed, and an explicit solution to the optimal controller is given. Due to the dependence of $\theta(k)$ on its past values, the new version of the maximum principle for the LQR problem needs to be established which can viewed as a generalization to the result for multiplicative noise systems \cite{ZWL12}.

\begin{lemma}
According to the linear system (\ref{ff1}) and the performance index (\ref{f2}). If the LQR problem $\min J_N$ is solvable, then the optimal ${\cal{G}}_{k}$-measurable control $u(k)$ satisfies the following equation
\begin{eqnarray}
0=\mbox{E}[B_{\theta(k)}'\eta_k+R_{\theta(k)}u(k)|{\cal{G}}_{k}],k=0,\cdots,N,\label{f3}
\end{eqnarray}
where the costate $\eta_k$ satisfies the following equation
\begin{eqnarray}
\eta_N&=&\mbox{E}[P_{\theta(N+1)}x(N+1)|{\cal{G}}_N],\label{f4}\\
\eta_{k-1}&=&\mbox{E}[A_{\theta(k)}'\eta_k+Q_{\theta(k)}x(k)|{\cal{G}}_{k-1}],k=0,\cdots,N.\label{f5}
\end{eqnarray}
\end{lemma}

\emph{Proof}. Denote $N$ as the final control horizon. It is known that $u(k)$ is ${\cal{G}}_{k}$-measurable. Consider the increment of the control variable $u(k)$ and deduce an expression of the corresponding variation of the performance index (\ref{f2})
\begin{eqnarray}
dJ_N&=&\mbox{E}[2x(N+1)'P_{\theta(N+1)}dx(N+1)+2\sum_{k=0}^Nx(k)'Q_{\theta(k)}dx(k)\nonumber\\
&&+2\sum_{k=0}^Nu(k)'R_{\theta(k)}du(k)]\nonumber\\
&=&\mbox{E}\{2x(N+1)'P_{\theta(N+1)}[F_x(N,0)dx_0+\sum_{i=0}^NF_x(N,i+1)B_{\theta(i)}du(i)]\nonumber\\
&&+2\sum_{k=d}^Nu(k-d)'R_{\theta(k)}du(k)\nonumber\\
&&+2\sum_{k=0}^Nx(k)'Q_{\theta(k)}[F_x(k-1,0)dx_0+\sum_{i=0}^{k-1}F_x(k-1,i+1)B_{\theta(i)}du(i)]\}\nonumber\\
&=&\mbox{E}\{2x(N+1)'P_{\theta(N+1)}[F_x(N,0)dx_0+\sum_{i=0}^NF_x(N,i+1)B_{\theta(i)}du(i)]\nonumber\\
&&+2\sum_{i=0}^Nu(i)'R_{\theta(i)}du(i)+2\sum_{k=0}^Nx(k)'Q_{\theta(k)}F_x(k-1,0)dx_0\nonumber\\
&&+2\sum_{i=0}^{N-1}\sum_{k={i+1}}^Nx(k)'Q_{\theta(k)}F_x(k-1,i+1)B_{\theta(i)}du(i)\},\label{f11}
\end{eqnarray}
where
\begin{eqnarray}
F_x(k,i)&=&A_{\theta(k)}\cdots A_{\theta(i)},i=0,\cdots,k,\nonumber\\
F_x(k,k+1)&=&I.\label{f8}
\end{eqnarray}
%In view of system (\ref{ff1}), we have
%\begin{eqnarray}
%dx(k+1)=F_x(k,0)dx_0+\sum_{i=0}^kF_x(k,i+1)B_{\theta(i)}du(i),\label{f7}
%\end{eqnarray}
%where
%\begin{eqnarray}
%F_x(k,i)&=&A_{\theta(k)}\cdots A_{\theta(i)},i=0,\cdots,k,\nonumber\\
%F_x(k,k+1)&=&I.\label{f8}
%\end{eqnarray}
%Plugging the equation (\ref{f7}) in (\ref{f6}) we deduce that
%\begin{eqnarray}
%dJ_N&=&\mbox{E}\{2x(N+1)'P_{\theta(N+1)}[F_x(N,0)dx_0+\sum_{i=0}^NF_x(N,i+1)B_{\theta(i)}du(i)]\nonumber\\
%&&+2\sum_{k=d}^Nu(k-d)'R_{\theta(k)}du(k)\nonumber\\
%&&+2\sum_{k=0}^Nx(k)'Q_{\theta(k)}[F_x(k-1,0)dx_0+\sum_{i=0}^{k-1}F_x(k-1,i+1)B_{\theta(i)}du(i)]\}.\label{f9}
%\end{eqnarray}
%It is easy to check that the last term of (\ref{f9}) can be rewritten as
%\begin{eqnarray}
%&&\sum_{k=0}^Nx(k)'Q_{\theta(k)}\sum_{i=0}^{k-1}F_x(k-1,i+1)B_{\theta(i)}du(i)\nonumber\\
%&=&\sum_{k=0}^N\sum_{i=0}^{k-1}x(k)'Q_{\theta(k)}F_x(k-1,i+1)B_{\theta(i)}du(i)\nonumber\\
%&=&\sum_{i=0}^{N-1}\sum_{k=i+1}^Nx(k)'Q_{\theta(k)}F_x(k-1,i+1)B_{\theta(i)}du(i).\label{f10}
%\end{eqnarray}
%Using (\ref{f10}), we deduce that
%\begin{eqnarray}
%dJ_N&=&\mbox{E}\{2x(N+1)'P_{\theta(N+1)}[F_x(N,0)dx_0+\sum_{i=0}^NF_x(N,i+1)B_{\theta(i)}du(i)]\nonumber\\
%&&+2\sum_{i=0}^Nu(i)'R_{\theta(i)}du(i)+2\sum_{k=0}^Nx(k)'Q_{\theta(k)}F_x(k-1,0)dx_0\nonumber\\
%&&+2\sum_{i=0}^{N-1}\sum_{k={i+1}}^Nx(k)'Q_{\theta(k)}F_x(k-1,i+1)B_{\theta(i)}du(i)\}.\label{f11}
%\end{eqnarray}
Since we just pay attention to the increment of $J_N$ caused by the increment of $u(i)$, the initial state $x_0$ is fixed and its increment $dx_0$ is thus $0$. Therefore,
\begin{eqnarray}
dJ_N
%&=&\mbox{E}\{2x(N+1)'P_{\theta(N+1)}\sum_{i=0}^NF_x(N,i+1)B_{\theta(i)}du(i)+2\sum_{i=0}^Nu(i)'R_{\theta(i)}du(i)\nonumber\\
%&&+2\sum_{i=0}^{N-1}\sum_{k=i+1}^Nx(k)'Q_{\theta(k)}F_x(k-1,i+1)B_{\theta(i)}du(i)\}\nonumber\\
&=&\mbox{E}\{2[x(N+1)'P_{\theta(N+1)}F_x(N,N+1)B_{\theta(N)}+u(N)'R_{\theta(N)}]du(N)\nonumber\\
&&+2\sum_{i=0}^{N-1}[x(N+1)'P_{\theta(N+1)}F_x(N,i+1)B_{\theta(i)}+u(i)'R_{\theta(i)}\nonumber\\
&&+\sum_{k=i+1}^Nx(k)'Q_{\theta(k)}F_x(k-1,i+1)B_{\theta(i)}du(i)\}.\label{f12}
\end{eqnarray}
Define
\begin{eqnarray}
\eta_i=\mbox{E}\{\sum_{k=i+1}^NF_x'(k-1,i+1)Q_{\theta(k)}x(k)+F_x'(N,i+1)P_{\theta(N+1)}x(N+1)|{\cal{G}}_i\},\label{f13}
\end{eqnarray}
then we have
\begin{eqnarray}
\eta_{i-1}&=&\mbox{E}\{\sum_{k=i}^NF_x'(k-1,i)Q_{\theta(k)}x(k)+F_x'(N,i)P_{\theta(N+1)}x(N+1)|{\cal{G}}_{i-1}\}\nonumber\\
%&=&\mbox{E}\{F_x'(i-1,i)Q_{\theta(i)}x(i)+\sum_{k=i+1}^NF_x'(k-1,i)Q_{\theta(k)}x(k)\nonumber\\
%&&+F_x'(N,i)P_{\theta(N+1)}x(N+1)|{\cal{G}}_{i-1}\}\nonumber\\
%&=&\mbox{E}\{Q_{\theta(i)}x(i)+\sum_{k=i+1}^NA_{\theta(i)}'F_x'(k-1,i+1)Q_{\theta(k)}x(k)\nonumber\\
%&&+A_{\theta(i)}'F_x'(N,i+1)P_{\theta(N+1)}x(N+1)|{\cal{G}}_{i-1}\}\nonumber\\
%&=&\mbox{E}\{Q_{\theta(i)}x(i)+A_{\theta(i)}'\mbox{E}[\sum_{k=i+1}^NF_x'(k-1,i+1)Q_{\theta(k)}x(k)\nonumber\\
%&&+F_x'(N,i+1)P_{\theta(N+1)}x(N+1)|{\cal{G}}_{i}]|{\cal{G}}_{i-1}\}\nonumber\\
&=&\mbox{E}\{Q_{\theta(i)}x(i)+A_{\theta(i)}'\eta_i|{\cal{G}}_{i-1}\}\nonumber.
\end{eqnarray}
%It has been shown that
%\begin{eqnarray}
%\lambda_{i-1}&=&\mbox{E}\{Q_{\theta(i)}x(i)+A_{\theta(i)}'\lambda_i|{\cal{G}}_{i-1}\},\label{f14}\\
%\lambda_{N}&=&\mbox{E}\{P_{\theta(N+1)}x(N+1)|{\cal{G}}_N\}.\label{f15}
%\end{eqnarray}
Based on (\ref{f13}), we deduce that
\begin{eqnarray}
dJ_N
%&=&\mbox{E}\{2\sum_{i=0}^N\mbox{E}[u(i)'R_{\theta(i)}du(i)|{\cal{G}}_i]\nonumber\\
%&&+2\mbox{E}[x(N+1)'P_{\theta(N+1)F_x(N,N+1)}|{\cal{G}}_N]B_{\theta(N)}du(N)\nonumber\\
%&&+2\sum_{i=0}^{N_1}\mbox{E}[x(N+1)P_{\theta(N+1)}F_x(N,i+1)\nonumber\\
%&&+\sum_{k=i+1}^Nx(k)'Q_{\theta(k)}F_x(k-1,i+1)|{\cal{G}}_i]B_{\theta(i)}du(i)\}\nonumber\\
%&=&\mbox{E}\{2\sum_{i=0}^Nu(i)'R_{\theta(i)}du(i)+2\sum_{i=0}^N\lambda_i'B_{\theta(i)}du(i)\}\nonumber\\
%&=&\mbox{E}\{2\sum_{i=0}^Nu(i)'R_{\theta(i)}du(i)+2\sum_{i=d}^N\lambda_i'B_{\theta(i)}du(i)\}\nonumber\\
%&=&\mbox{E}\{2\sum_{i=0}^N[\lambda_i'B_{\theta(i)}+u(i)'R_{\theta(i)}]du(i)\}\nonumber\\
&=&\mbox{E}\{2\sum_{i=0}^N\mbox{E}[\eta_i'B_{\theta(i)}+u(i)'R_{\theta(i)}|{\cal{G}}_{i}]du(i)\}.\label{f16}
\end{eqnarray}
It concludes from (\ref{f16}) that the necessary condition for the minimum can be given as follows
\begin{eqnarray}
\mbox{E}[\eta_i'B_{\theta(i)}+u(i)'R_{\theta(i)}|{\cal{G}}_{i}]=0,i=0,\cdots,N.\label{f17}
\end{eqnarray}
This completes the proof.

%%%%%%%%%%%%%%%%%%%%%%%%%%%%%%%%%%%%%%%%%%%%%%%%%%%%%%%%%%%%%%%%%%%%%%%%%%%%%%%%%%%%%%%%%%%%%%%%%%%%%%%%%%%%%%%%%%%%%%%%%%%%%%%%%%%%%%%%%%%%%%%%%%%%%%%%%%%%%%%%
%\subsection{Solution to the finite-horizon LQR}

In what follows, we will derive the analytic solution for the LQR problem, and give the necessary and sufficient conditions for the existence of the optimal controller.

\begin{theorem}
Problem 1 has a unique solution if and only if the following coupled difference equations
\begin{eqnarray}
\Upsilon_{i}(k)&=&B_{i}'(\sum_{{j}=1}^L\lambda_{i,j}P_{j}(k+1))B_{i}+R_{i},\label{f18}\\
M_{i}(k)&=&B_{i}'(\sum_{j=1}^L\lambda_{i,j}P_{j}(k+1))A_{i},\label{f19}\\
P_{i}(k)&=&A_{i}'(\sum_{j=1}^L\lambda_{i,j}P_{j}(k+1))A_{i}+Q_{i}-M_{i}(k)'\Upsilon_{i}(k)^{-1}M_{i}(k),\label{f20}
%P_{l_N}(N+1)&=&\sum_{l_{N+1}=1}^L\lambda_{{l_N},{l_{N+1}}}P_{l_{N+1}}(N+1)\nonumber
\end{eqnarray}
are well defined for $k=N,\cdots,0, i=1,\cdots,L$, that is $\Upsilon_{i}(k), k=N,\cdots,0, i=1,\cdots,L$ are all invertible. If this condition is satisfied, the analytical solution to the optimal control can be given as
\begin{eqnarray}
u(k)=-\Upsilon_{i}(k)^{-1}M_{i}(k){x}(k), i=1,\cdots,L,\label{f21}
\end{eqnarray}
for $k=N,\cdots,0$.
%\begin{eqnarray}
%\hat{x}(k|k-1)=\mbox{E}[x(k)|{\cal{G}}_{k-1}]=A_{\theta(k-1)}x(k-1)+B_{\theta(k-1)}u(k-2).\label{f22}
%\end{eqnarray}
The corresponding optimal performance index is given by
\begin{eqnarray}
J_N=\mbox{E}[x(0)'P_{\theta(0)}(0)x(0)].\label{f23}
\end{eqnarray}
The solution to FBMJDE (\ref{f3})-(\ref{f5}), i.e., the relationship of $\eta_{k-1}$ and $x(k)$, is given as
\begin{eqnarray}
\eta_{k-1}=(\sum_{j=1}^L\lambda_{i,j}P_{j}(k))x(k), i=1,\cdots,L.\label{f24}
\end{eqnarray}
\end{theorem}

\emph{Proof}. ``Necessary": Assume that Problem 1 has a unique solution. By the induction, we will prove that $\Upsilon_{i}(k)$ in (\ref{f18}) is invertible for all $k=N,\cdots,0, i=1,\cdots,L$, and $u(k)$ satisfies (\ref{f21}). Define
\begin{eqnarray}
J(k)\stackrel{\triangle}{=} \mbox{E}\{\sum_{i=k}^N(x(i)'Q_{\theta(i)}x(i)+u(i)'R_{\theta(i)}u(i))+x(N+1)'P_{\theta(N+1)}(N+1)x(N+1)|{\cal{G}}_{k}\},\label{f25}
\end{eqnarray}
for $k=N,\cdots,0$. For $k=N$, (\ref{f25}) becomes
\begin{eqnarray}
J(N){=} \mbox{E}\{x(N)'Q_{\theta(N)}x(N)+u(N)'R_{\theta(N)}u(N)+x(N+1)'P_{\theta(N+1)}(N+1)x(N+1)|{\cal{G}}_{N}\},\label{f26}
\end{eqnarray}
Based on (\ref{ff1}), we deduce that $J(N)$ can be represented as quadratic function of $x(N)$ and $u(N)$. The uniqueness of the optimal controller $u(N)$ indicates that the quadratic term of $u(N)$ is positive for any nonzero $u(N)$. Let $x(N)=0$ and substitute (\ref{ff1}) in (\ref{f26}), we have
\begin{eqnarray}
J(N)&=&\mbox{E}\{u(N)'(R_{\theta(N)}+B_{\theta(N)}'P_{\theta(N+1)}(N+1)B_{\theta(N)})u(N)|{\cal{G}}_{N}\}\nonumber\\
%&=&\mbox{E}\{\mbox{E}\{u(N)'(R_{\theta(N)}+B_{\theta(N)}'P_{\theta(N+1)}B_{\theta(N)})u(N)|{\cal{G}}_{N}\}|{\cal{G}}_{N}\}\nonumber\\
&=&u(N)'[R_{i}+B_{i}'(\sum_{j=1}^L\lambda_{i,j}P_{j}(N+1))B_{i}]u(N)\nonumber\\
%&=&u(N-1)'[\sum_{l_N=1}^L\lambda_{l_{N-1},l_N}(R_{l_N}+B_{l_N}'P_{l_N}(N+1)B_{l_N})]u(N-1)\nonumber\\
&=&u(N)'\Upsilon_{i}(N)u(N)>0,i=1,\cdots,L,\label{f27}
\end{eqnarray}
where
\begin{eqnarray}
\Upsilon_{i}(N)=R_{i}+B_{i}'(\sum_{j=1}^L\lambda_{i,j}P_{j}(N+1))B_{i}. \label{f28}
\end{eqnarray}
It can be concluded that $\Upsilon_{i}(N)>0$.

In what follows, the optimal controller $u(N)$ is to be calculated. Applying (\ref{ff1}), (\ref{f3}) and (\ref{f4}), we obtain that
\begin{eqnarray}
0&=&\mbox{E}[B_{\theta(N)}'\eta_N+R_{\theta(N)}u(N)|{\cal{G}}_{N}]\nonumber\\
&=&\mbox{E}[B_{\theta(N)}'\mbox{E}(P_{\theta(N+1)}(N+1)x(N+1)|{\cal{G}}_N)+R_{\theta(N)}u(N)|{\cal{G}}_{N}]\nonumber\\
&=&\mbox{E}[B_{\theta(N)}'\mbox{E}[P_{\theta(N+1)}(N+1)(A_{\theta(N)}x(N)+B_{\theta(N)}u(N))|{\cal{G}}_N]+R_{\theta(N)}u(N)|{\cal{G}}_{N}]\nonumber\\
&=&B_{i}'(\sum_{j=1}^L\lambda_{i,j}P_{j}(N+1))A_{i}x(N)+B_{i}'(\sum_{j=1}^L\lambda_{i,j}P_{j}(N+1))B_{i}u(N)\nonumber\\
&&+R_{i}u(N)\nonumber\\
&=&B_{i}'(\sum_{j=1}^L\lambda_{i,j}P_{j}(N+1))A_{i}x(N)+[B_{i}'(\sum_{j=1}^L\lambda_{i,j}P_{j}(N+1))B_{i}\nonumber\\
&&+R_{i}]u(N).\nonumber
\end{eqnarray}
It follows from the above equation that
\begin{eqnarray}
u(N)=-\Upsilon_{i}(N)^{-1}M_{i}(N){x}(N),i=1,\cdots,L,\label{f29}
\end{eqnarray}
where $\Upsilon_{i}(N)$ is as in (\ref{f28}) and $M_{i}(N)$ is as follows
\begin{eqnarray}
M_{i}(N)=B_{i}'(\sum_{j=1}^L\lambda_{i,j}P_{j}(N+1))A_{i}.\label{f30}
\end{eqnarray}

In the following, we will show that $\eta_{N-1}$ is with the form as (\ref{f24}). In view of (\ref{ff1}), (\ref{f5}), and (\ref{f29}), one yields
\begin{eqnarray}
\eta_{N-1}&=&\mbox{E}\{A_{\theta(N)}'\eta_N+Q_{\theta(N)}x(N)|{\cal{G}}_{N-1}\}\nonumber\\
&=&\mbox{E}\{A_{\theta(N)}'\mbox{E}[P_{\theta(N+1)}(N+1)x(N+1)|{\cal{G}}_N]+Q_{\theta(N)}x(N)|{\cal{G}}_{N-1}\}\nonumber\\
&=&\mbox{E}\{A_{\theta(N)}'\mbox{E}[P_{\theta(N+1)}(N+1)(A_{\theta(N)}x(N)+B_{\theta(N)}u(N))|{\cal{G}}_N]+Q_{\theta(N)}x(N)|{\cal{G}}_{N-1}\}\nonumber\\
%&=&\mbox{E}\{A_{l_N}'P_{l_N}(N+1)A_{l_N}x(N)+A_{l_N}'P_{l_N}(N+1)B_{l_N}u(N)+Q_{l_N}x(N)|{\cal{G}}_{N-1}\}\nonumber\\
&=&\mbox{E}\{[A_{i}'(\sum_{j=1}^L\lambda_{i,j}P_{j}(N+1))A_{i}+Q_{i}]{x}(N)+A_{i}'(\sum_{j=1}^L\lambda_{i,j}P_{j}(N+1))B_{i}u(N)|{\cal{G}}_{N-1}\}\nonumber\\
&=&\mbox{E}\{[A_{i}'(\sum_{j=1}^L\lambda_{i,j}P_{j}(N+1))A_{i}+Q_{i}]{x}(N)-M_{i}(N)'\Upsilon_{i}(N)^{-1}M_{i}(N){x}(N)|{\cal{G}}_{N-1}\}\nonumber\\
&=&(\sum_{i=1}^L\lambda_{s,i}P_{i}(N))x(N),s=1,\cdots,L,\nonumber
\end{eqnarray}
where
\begin{eqnarray}
P_{i}(N)=A_{i}'(\sum_{j=1}^L\lambda_{i,j}P_{j}(N+1))A_{i}+Q_{i}-M_{i}(N)'\Upsilon_{i}(N)^{-1}M_{i}(N).\nonumber
\end{eqnarray}

To proceed the induction proof, we take any $n$ with $1\leq n\leq N$, and assume that $\Upsilon_{i}(k)(i=1,\cdots,L)$ is invertible and that the optimal controller $u(k)$ and the optimal costate $\eta_{k-1}$ are as (\ref{f21}) and (\ref{f24}) for all $k\geq n+1$. In the next, it needs to show that these conditions will also be satisfied for $k=n$. Follow the similar derivation procedure for $\Upsilon_{i}(N)(i=1,\cdots,L)$ and let $x(n)=0$, we will check the quadratic term of $u(n)$ in $J(n)$. In view of (\ref{ff1}), (\ref{f3}), and (\ref{f5}) for $k\geq n+1$, we have
\begin{eqnarray}
&&\mbox{E}\{x(k)'\eta_{k-1}-x(k+1)'\eta_k|{\cal{G}}_{n+1}\}\nonumber\\
&=&\mbox{E}\{x(k)'\mbox{E}[Q_{\theta(k)}x(k)+A_{\theta(k)}'\eta_k|{\cal{G}}_{k-1}]-[A_{\theta(k)}x(k)+B_{\theta(k)}u(k)]'\eta_k|{\cal{G}}_{n+1}\}\nonumber\\
&=&\mbox{E}\{\mbox{E}[x(k)'Q_{\theta(k)}x(k)+x(k)'A_{\theta(k)}'\eta_k|{\cal{G}}_{k-1}]-\mbox{E}[x(k)'A_{\theta(k)}'\eta_k
+u(k)'B_{\theta(k)}'\eta_k|{\cal{G}}_{k-1}]|{\cal{G}}_{n+1}\}\nonumber\\
&=&\mbox{E}\{\mbox{E}[x(k)'Q_{\theta(k)}x(k)-u(k)'B_{\theta(k)}'\eta_k|{\cal{G}}_{k-1}]|{\cal{G}}_{n+1}\}\nonumber\\
&=&\mbox{E}\{x(k)'Q_{\theta(k)}x(k)-u(k)'B_{\theta(k)}'\eta_k|{\cal{G}}_{n+1}\}\nonumber\\
&=&\mbox{E}\{x(k)'Q_{\theta(k)}x(k)+u(k)'R_{\theta(k)}u(k)|{\cal{G}}_{n+1}\}\nonumber
\end{eqnarray}
Adding from $k=n+1$ to $k=N$ on both sides of the above equation, we obtain that
\begin{eqnarray}
&&\mbox{E}\{x(n+1)'\eta_n-x(N+1)'\eta_N|{\cal{G}}_{n+1}\}\nonumber\\
&=&\sum_{k=n+1}^N\mbox{E}\{x(k)'\eta_{k-1}-x(k+1)'\eta_k|{\cal{G}}_{n+1}\}\nonumber\\
&=&\sum_{k=n+1}^N\mbox{E}\{x(k)'Q_{\theta(k)}x(k)+u(k)'R_{\theta(k)}u(k)|{\cal{G}}_{n+1}\}.\label{f31}
\end{eqnarray}
So we have from (\ref{f31}) that
\begin{eqnarray}
J(n)&=&\mbox{E}\{\sum_{i=n}^N(x(i)'Q_{\theta(i)}x(i)+u(i)'R_{\theta(i)}u(i))\nonumber\\
&&+x(N+1)'P_{\theta(N+1)}x(N+1)|{\cal{G}}_{n}\}\nonumber\\
&=&\mbox{E}\{x(n)'Q_{\theta(n)}x(n)+u(n)'R_{\theta(n)}u(n)\nonumber\\
&&+\mbox{E}[\sum_{i=n+1}^N(x(i)'Q_{\theta(i)}x(i)+u(i)'R_{\theta(i)}u(i))\nonumber\\
&&+x(N+1)'P_{\theta(N+1)}x(N+1)|{\cal{G}}_{n+1}]|{\cal{G}}_{n}\}\nonumber\\
&=&\mbox{E}\{x(n)'Q_{\theta(n)}x(n)+u(n)'R_{\theta(n)}u(n)+\mbox{E}[x(n+1)'\eta_n\nonumber\\
&&-x(N+1)'\eta_N+x(N+1)'P_{\theta(N+1)}x(N+1)|{\cal{G}}_{n+1}]|{\cal{G}}_{n}\}\nonumber\\
&=&\mbox{E}\{x(n)'Q_{\theta(n)}x(n)+u(n)'R_{\theta(n)}u(n)+\mbox{E}[x(n+1)'\eta_n|{\cal{G}}_{n+1}]|{\cal{G}}_{n}\}\nonumber\\
&=&\mbox{E}\{x(n)'Q_{\theta(n)}x(n)+u(n)'R_{\theta(n)}u(n)+x(n)'A_{\theta(n)}'\eta_n\nonumber\\
&&+u(n)'B_{\theta(n)}'\eta_n|{\cal{G}}_{n}\}\nonumber\\
&=&\mbox{E}\{u(n)'R_{\theta(n)}u(n)+u(n)'B_{\theta(n)}'\eta_n|{\cal{G}}_{n}\}.\label{f32}
\end{eqnarray}
Note that
\begin{eqnarray}
\eta_n&=&(\sum_{j=1}^L\lambda_{i,j}P_{j}(n+1))x(n+1)\nonumber\\
&=&(\sum_{j=1}^L\lambda_{i,j}P_{j}(n+1))[A_{i}x(n)+B_{i}u(n)]. \label{f33}
\end{eqnarray}
Substitute (\ref{f33}) in (\ref{f32}), we deduce that
\begin{eqnarray}
J(n)&=&\mbox{E}\{u(n)'R_{i}u(n)+u(n)'B_{i}'(\sum_{j=1}^L\lambda_{i,j}P_{j}(n+1))B_{i}u(n)|{\cal{G}}_{n}\}\nonumber\\
&=&u(n)'[B_{i}'(\sum_{j=1}^L\lambda_{i,j}P_{j}(n+1))B_{i}+R_{i}]u(n)\nonumber\\
&=&u(n)'\Upsilon_{i}(n)u(n).\label{f34}
\end{eqnarray}
It is concluded from the uniqueness of the optimal controller that $J(n)$ must be positive for any $u(n)\neq 0$. So we have $\Upsilon_{i}(n)>0, i=1,\cdots,L$.

To derive the optimal controller $u(n)$, plugging (\ref{f33}) in (\ref{f3}) yields
\begin{eqnarray}
0&=&\mbox{E}\{B_{i}'\eta_n+R_{i}u(n)|{\cal{G}}_{n}\}\nonumber\\
&=&B_{i}'(\sum_{j=1}^L\lambda_{i,j}P_{j}(n+1))A_{i}x(n)\nonumber\\
&&+(B_{i}'(\sum_{j=1}^L\lambda_{i,j}P_{j}(n+1))B_{i}+R_{i})u(n)\nonumber\\
%&=&\mbox{E}\{B_{l_n}'(\sum_{l_{n+1}=1}^L\lambda_{l_n,l_{n+1}}P_{l_{n+1}}(n+1))A_{l_n}x(n)\nonumber\\
%&&-(B_{l_n}'(\sum_{l_{n+1}=1}^L\lambda_{l_n,l_{n+1}}P_{l_{n+1}}(n+1))B_{l_n}+R_{l_n})u(n)|{\cal{G}}_{n}\}\nonumber\\
%&=&\sum_{l_n=1}^L\lambda_{l_{n-1},l_n}(B_{l_n}'P_{l_n}(n+1)A_{l_n})\hat{x}(n|n-1)+\sum_{l_n=1}^L\lambda_{l_{n-1},l_n}(B_{l_n}'P_{l_n}(n+1)B_{l_n}+R_{l_n})u(n-1)\nonumber\\
&=&M_{i}(n){x}(n)+\Upsilon_{i}(n)u(n).\nonumber
\end{eqnarray}
Using the above equation, we get
\begin{eqnarray}
u(n)=-\Upsilon_{i}(n)^{-1}M_{i}(n){x}(n),\label{f35}
\end{eqnarray}
where
\begin{eqnarray}
\Upsilon_{i}(n)&=&B_{i}'(\sum_{j=1}^L\lambda_{i,j}P_{j}(n+1))B_{i}+R_{i},\nonumber\\
M_{i}(n)&=&B_{i}'(\sum_{j=1}^L\lambda_{i,j}P_{j}(n+1))A_{i}.\nonumber
\end{eqnarray}
Now, we proceed to derive that $\eta_{n-1}$ is of the form as (\ref{f24}). In terms of (\ref{f5}), (\ref{f33}) and (\ref{f35}), we have
\begin{eqnarray}
\eta_{n-1}&=&\mbox{E}\{Q_{\theta(n)}x(n)+A_{\theta(n)}'\eta_n|{\cal{G}}_{n-1}\}\nonumber\\
&=&\mbox{E}\{Q_{\theta(n)}x(n)+A_{\theta(n)}'(\sum_{j=1}^L\lambda_{\theta(n),j}P_{j}(n+1))(A_{\theta(n)}x(n)+B_{\theta(n)}u(n))|{\cal{G}}_{n-1}\}\nonumber\\
&=&\sum_{i=1}^L\lambda_{s,i}(Q_{i}+A_{i}'(\sum_{j=1}^L\lambda_{i,j}P_{j}(n+1))A_{i}){x}(n)\nonumber\\
&&+\sum_{i=1}^L\lambda_{s,i}A_{i}'(\sum_{j=1}^L\lambda_{i,j}P_{j}(n+1))B_{i}u(n)\nonumber\\
&=&\sum_{i=1}^L\lambda_{s,i}\{(Q_{i}+A_{i}'(\sum_{j=1}^L\lambda_{i,j}P_{j}(n+1))A_{i})\nonumber\\
&&-M_{i}(n)'\Upsilon_{i}(n)^{-1}M_{i}(n)\}{x}(n)\nonumber\\
&=&\{\sum_{i=1}^L\lambda_{s,i}P_{i}(n)\}{x}(n),\nonumber
\end{eqnarray}
where
\begin{eqnarray}
P_{i}(n)&=&Q_{i}+A_{i}'(\sum_{j=1}^L\lambda_{i,j}P_{j}(n+1))A_{i}-M_{i}(n)'\Upsilon_{i}(n)^{-1}M_{i}(n).\nonumber
\end{eqnarray}

``Sufficiency": Suppose $\Upsilon_{i}(k)>0,i=1,\cdots,L$, then we will prove that Problem 1 a has a unique solution. Define
\begin{eqnarray}
V_N(k,x(k))\stackrel{\triangle}{=} \mbox{E}[x(k)'(\sum_{i=1}^L\lambda_{s,i}P_{i}(k)){x}(k)]=\mbox{E}[x(k)'P_{\theta(k)}(k){x}(k)].\label{f36}
\end{eqnarray}
Applying (\ref{f36}), (\ref{f18})-(\ref{f20}), we deduce that
\begin{eqnarray}
&&V_N(k,x(k))-V_N(k+1,x(k+1))\nonumber\\
&=&\mbox{E}\{x(k)'Q_{i}x(k)+u(k)'R_{i}u(k)-[u(k)+\Upsilon_{i}(k)^{-1}M_{i}(k){x}(k)]\nonumber\\
&&\times \Upsilon_{i}(k)[u(k)+\Upsilon_{i}(k)^{-1}M_{i}(k){x}(k)]\}.\label{f37}
\end{eqnarray}
Adding from $k=0$ to $k=N$ on both sides of (\ref{f37}), the performance index (\ref{f2}) is rewritten as
\begin{eqnarray}
J_N&=&\mbox{E}\{x(0)'P_{\theta(0)}(0)x(0)+\sum_{k=0}^N[u(k)+\Upsilon_{i}(k)^{-1}M_{i}(k){x}(k)]\nonumber\\
&&\times \Upsilon_{i}(k)[u(k)+\Upsilon_{i}(k)^{-1}M_{i}(k){x}(k)]\}.\nonumber
\end{eqnarray}

Note that $\Upsilon_{i}(k)>0, i=1,\cdots,L$. Thus Problem 1 has a unique solution, and the optimal controller is given by
\begin{eqnarray}
u(k)=-\Upsilon_{i}(k)^{-1}M_{i}(k){x}(k).\nonumber
\end{eqnarray}
The corresponding optimal performance index is given by
\begin{eqnarray}
J_N=\mbox{E}[x(0)'P_{\theta(0)}(0)x(0)].\nonumber
\end{eqnarray}
This completes the proof.

\begin{remark}
Necessary and sufficient condition for the existence of the discrete-time Markov jump linear LQR problem is given in Theorem 1, in which it just requires the input penalty matrix $R$ is positive semi-definite. It can be found that the existed results usually consider the case that $R$ is positive definite \cite{BS75}, \cite{CFM05}, \cite{CABM99}, \cite{CT03} or a set of matrix expressions is positive definite \cite{CWC86}, \cite{JC88}, and only sufficient conditions for the existence of the LQR controller is given. In this paper, an analytical solution to the forward-backward
Markov jumping parameter difference equation associated with
optimal control is presented. This forms the basis on which we
supply the necessary and sufficient conditions for the existence the optimal LQR controller for MJLS.
\end{remark}

%\begin{remark}
%Compared to the result for discrete-time multiplicative noise systems, the result developed in Theorem 1 becomes more general since the jumping parameters are correlated at same and adjoining time. If the the jumping parameter is independent between the adjacent time, the result developed in Theorem 1 is reduced to the result developed in [18] for the multiplicative noise system without delays. So the result developed in Theorem 1 can be considered as the generalization for the multiplicative case.
%\end{remark}

%%%%%%%%%%%%%%%%%%%%%%%%%%%%%%%%%%%%%%%%%%%%%%%%%%%%%%%%%%%%%%%%%%%%%%%%%%%%%%%%%%%%%%%%%%%%%%%%%%%%%%%%%%%%%%%%%%%%%%%%%%%%%%%%%%%%%%%%%%%%%%%%%%%%%%%%%%%%%%
\section{Infinite-Horizon LQR for MJLS}
\subsection{Problem statement}
For the infinite horizon quadratic optimal control problems to be analyzed in this section, we consider a time-invariant version of the model (\ref{ff1}). We will be interested in the problem of minimizing the infinite horizon cost function given by
\begin{eqnarray}
J&=&\mbox{E}\{\sum_{k=0}^{\infty}[x(k)'Q_{\theta(k)}x(k)+u(k)'R_{\theta(k)}u(k)]\}.\label{f38}
\end{eqnarray}
\begin{definition}
We say that the linear system with Markov jump parameter (\ref{f1}) with $u(k)=0$ is mean square stable (MSS) if for any initial condition $x_0$ and $\theta(0)$, there holds
\begin{eqnarray}
\lim_{k\rightarrow \infty}\mbox{E}(x(k)'x(k))=0.\nonumber
\end{eqnarray}
\end{definition}
\begin{definition}
We say that system (\ref{ff1}) is mean square stabilizable if there is a ${\cal{G}}_k$-measurable controller $u(k)=F_{\theta(k)}x(k)$ satisfying $\lim_{k\rightarrow \infty}\mbox{E}[u(k)'u(k)]=0$, such that system (\ref{ff1}) is asymptotically mean square stable.
\end{definition}
\begin{definition}
The following MJLS
\begin{eqnarray}
x(k+1)&=&A_{\theta(k)}x(k),y(k)=C_{\theta(k)}x(k)\label{f39}
\end{eqnarray}
is said to be exactly observable, if for any $N$
\begin{eqnarray}
y(k)=0,a.s. \forall s\leq k\leq N\Rightarrow x_0=0.\nonumber
\end{eqnarray}
\end{definition}
Denote $A=(A_1,\cdots,A_L), B=(B_1,\cdots,B_L)$, and $C=(C_1,\cdots,C_L)$. For brevity, we usually say that the pair $(A, B)$ is mean square stabilizable if system (\ref{ff1}) is mean square stabilizable, and say that the pair $(C, A)$ is exactly observable if system (\ref{f39}) is exactly observable.

\emph{Problem 2}: Find the ${\cal{G}}_k$-measurable controller $u(k)=F_{\theta(k)}x(k),k\geq 0$, such that the closed loop system is asymptotically stable in the mean square sense, and the corresponding cost function (\ref{f38}) is minimized.

\begin{assumption}
$R_i (i=1,\cdots,L)$ is positive definite; $Q_i (i=1,\cdots,L)$ is positive semi-definite, that is, $Q_i=C_iC_i'(i=1,\cdots,L)$ for some matrix $C_i (i=1,\cdots,L)$.
\end{assumption}
\begin{assumption}
$(C,A)$ is exactly observable.
\end{assumption}
For clarity, we rewrite $\Upsilon_{i}(k), P_i(k)$, and $M_i(k) (i=1,\cdots,L)$ in (\ref{f18})-(\ref{f20}) as $\Upsilon_{i}^N(k), P_i^N(k)$, and $M_i^N(k) (i=1,\cdots,L)$. Without loss of generality, we set the terminal weight matrix $P_{j}^N(N+1) (j=1,\cdots,L)$ in the cost function to be zero.
Define the following coupled algebraic Riccati equation
\begin{eqnarray}
P_i&=&A_i'(\sum_{j=1}^L\lambda_{ij}P_j)A_i+Q_i-A_i'(\sum_{j=1}^L\lambda_{ij}P_j)B_i [B_i'(\sum_{j=1}^L\lambda_{ij}P_j)B_i+R_i]^{-1}\nonumber\\
&&\times B_i'(\sum_{j=1}^L\lambda_{ij}P_j)A_i,i=1,2,\cdots,L, \label{f40}\\
\Upsilon_i&=&B_i(\sum_{j=1}^L\lambda_{ij}P_j)B_i+R_i, i=1,2,\cdots,L, \label{f41}\\
M_i&=&B_i(\sum_{j=1}^L\lambda_{ij}P_j)A_i, i=1,2,\cdots,L.\label{f42}
\end{eqnarray}

\begin{lemma}
For any $N\geq 0$, $P_{i}^{N}(k)\geq 0$.
\end{lemma}
{ \emph{Proof}:} \ \  From (\ref{f41}) and (\ref{f42}), we have
\begin{eqnarray*}
[M_{i}^{N}(k)]'[\Upsilon_{i}^{N}(k)]^{-1}M_{i}^{N}(k)=-[M_{i}^{N}(k)]'K_{i}^{N}(k)-[K_{i}^{N}(k)]'M_{i}^{N}(k)
-[K_{i}^{N}(k)]'\Upsilon_{i}^{N}(k)K_{i}^{N}(k)
\end{eqnarray*}
in which $K_{i}^{N}(k)=-[\Upsilon_{i}^{N}(k)]^{-1}M_{i}^{N}(k)$.

In considering of (\ref{f40}), it yields that
\begin{eqnarray*}
P_{i}^{N}(k)&=&A_{i}'(\sum_{j=1}^L\lambda_{i,j}P_{j}^{N}(k+1))A_{i}
+Q_{i}+[M_{i}^{N}(k)]'K_{i}^{N}(k)\nonumber\\
&&+[K_{i}^{N}(k)]'M_{i}^{N}(k)
+[K_{i}^{N}(k)]'\Upsilon_{i}^{N}(k)K_{i}^{N}(k)\nonumber\\
&=&Q_{i}+[K_{i}^{N}(k)]'\Upsilon_{i}^{N}(k)K_{i}^{N}(k)+[A_{i}+B_{i}K_{i}^{N}(k)]'
(\sum_{j=1}^L\lambda_{i,j}P_{j}^{N}(k+1))[A_{i}+B_{i}K_{i}^{N}(k)].
\end{eqnarray*}
In view of the terminal condition $P_{j}^{N+1}(N+1)=0$ and $Q_{i}\geq0$, we can obtain that $P_{i}^{N}(N)\geq0$, and by induction, it is not hard to verify that $P_{i}^{N}(k)\geq0$, for $0\leq k \leq N$.

\begin{lemma}
When $R_{i}>0$, Problem 1 has a unique solution.
\end{lemma}
{ \emph{Proof}:} \ \  From Lemma 2, the expression of (\ref{f41}) and $R_{i}>0$, it is easy to obtain that $\Upsilon_{i}^{N}(k)>0$. According to the result of Theorem 1, in the case of $R_{i}>0$, we have Problem 1 has a unique solution.

\begin{theorem}
Under Assumptions 1 and 2, if the system (1) is mean square stabilizable , we have the following properties:

For any $k\geq 0$, $P_{i}^{N}(k)$ is convergent when $N\rightarrow\infty$, i.e.,
$\lim\limits_{N\rightarrow\infty}P_{i}^{N}(k)=P_{i}$, in which $P_{i}$ satisfies (\ref{f40})-(\ref{f42}), furthermore, $P_{i}>0$.
\end{theorem}

\emph{Proof}. First, we show that $P_{\theta(0)}^{N}(0)$ is increasing with respect to $N$. On the ground of $J_{N}\leq J_{N+1}$, we have that
$J_{N}^{\ast}\leq J_{N+1}^{\ast}$ for any initial value $x_{0}$. From (\ref{f23}), we obtain that
\begin{eqnarray*}
\mbox{E}[x(0)'P_{\theta(0)}^{N}(0)x(0)]\leq \mbox{E}[x(0)'P_{\theta(0)}^{N+1}(0)x(0)].
\end{eqnarray*}
In view of the arbitrary of $x(0)$, it implies that $P_{\theta(0)}^{N}(0)\leq P_{\theta(0)}^{N+1}(0)$, i.e., $P_{\theta(0)}^{N}(0)$ is increasing with respect to $N$.

Next, we will show the boundedness of $P_{\theta(0)}^{N}(0)$. When the system (\ref{ff1}) is stabilizable in the mean square sense, there exists $u(k)=F_{\theta(k)}x(k)$ satisfying
\begin{eqnarray*}
\lim\limits_{k\rightarrow\infty}E(x_{k}'x_{k})=0.
\end{eqnarray*}
Hence, we have that
\begin{eqnarray*}
J_{N}^{\ast}\leq J&=&\mbox{E}[\sum_{k=0}^\infty(x(k)'Q_{\theta(k)}x(k)+u(k)'R_{\theta(k)}u(k)')]\\
&=&\mbox{E}[\sum_{k=0}^\infty x(k)'Q_{\theta(k)}x(k)+x(k)'F_{\theta(k)}'R_{\theta(k)}F_{\theta(k)}x(k)]\\
&=&\mbox{E}[\sum_{k=0}^\infty x(k)'(Q_{\theta(k)}+F_{\theta(k)}'R_{\theta(k)}F_{\theta(k)})x(k)]\\
&&\leq \lambda\mbox{E}[\sum_{k=0}^\infty x(k)'x(k)]\\
&&\leq \lambda\cdot c \cdot \mbox{E}[x(0)'x(0)]
\end{eqnarray*}
in which $\lambda$ denotes the maximum eigenvalue  of $(Q_{\theta(k)}+F_{\theta(k)}'R_{\theta(k)}F_{\theta(k)})$ and $c$ is a positive constant. The above formula implies that
\begin{eqnarray*}
\mbox{E}[x(0)'P_{\theta(0)}^{N}(0)x(0)]\leq \lambda\cdot c \cdot \mbox{E}[x(0)'x(0)],
\end{eqnarray*}
i.e.,
\begin{eqnarray*}
P_{\theta(0)}^{N}(0)\leq \lambda \cdot c I.
\end{eqnarray*}
From now on, we can say that $P_{\theta(0)}^{N}(0)$ is bounded. In considering of the monotonicity of $P_{\theta(0)}^{N}(0)$, we deduce that $P_{\theta(0)}^{N}(0)$ is convergent. Note that the variables given in (\ref{f18})-(\ref{f20}) are time invariant for $N$ due to the choice that $P_{j}(N+1)=0,(j=1,\cdots,L)$, so we have
\begin{eqnarray*}
\lim\limits_{k\rightarrow\infty}P_{i}^{N}(k)=\lim\limits_{k\rightarrow\infty}P_{i}^{N-k}(0)=P_{i},i=1,\cdots,L.
\end{eqnarray*}
At the same time, we have that
\begin{eqnarray*}
\lim\limits_{k\rightarrow\infty}\Upsilon_{i}^{N}(k)=\Upsilon_{i},\ \ \lim\limits_{k\rightarrow\infty}M_{i}^{N}(k)=M_{i}.
\end{eqnarray*}

Now we will illustrate $P_{i}>0$. Since $J_{N}^{\ast}=\mbox{E}[x(0)'P_{\theta(0)}^{N}(0)x(0)]\geq0$, we can obtain that $P_{\theta(0)}^{N}(0)\geq 0$ for the arbitrary of $x(0)$. Next we mainly investigate that there exists a positive integer $N_0$ such that $P_{\theta_{(0)}}^{N_0}(0)>0$. If not, there must exist a nonempty set as follows:
\begin{eqnarray*}
Z_{N}=\left\{x\in R^{n}: x\neq 0, x'P_{\theta(0)}^{N}(0)x=0\right\}.
\end{eqnarray*}
The monotonically increasing of $P_{\theta(0)}^{N}(0)$ implies that if $\mbox{E}[x'P_{\theta(0)}^{N+1}(0)x]=0$, then $\mbox{E}[x'P_{\theta(0)}^{N}(0)x]=0$, i.e., $Z_{N+1}\subseteq Z_{N}$. As $Z_{N}$ is a nonempty finite dimensional set, thus
\begin{eqnarray*}
1\leq\cdots \leq dim(Z_{2})\leq dim(Z_{1})\leq dim(Z_{0})\leq n.
\end{eqnarray*}
Therefore, there exists a positive integer $N_{1}$ such that for any $N>N_{1}$ we have
\begin{eqnarray*}
dim(Z_{N})= dim(Z_{N_{1}}),
\end{eqnarray*}
i.e., $Z_{N}=Z_{N_{1}}$, furthermore, $\bigcap\limits_{N\geq0}Z_{N}=Z_{N_{1}}\neq0$.
Therefore, there exists a nonzero $x\in Z_{N}$ such that
$x'P_{\theta(0)}^{N}(0)x=0$.  Now let $x_{0}=x$, then
\begin{eqnarray*}
J_N^{\ast}&=&\mbox{E}[\sum_{k=0}^Nx^{\ast}(k)'Q_{\theta(k)}x^{\ast}(k)+u^{\ast}(k)'R_{\theta(k)}u^{\ast}(k)]\\
&=&\mbox{E}[x_{0}'P_{\theta(0)}^{N}(0)x_{0}]=0.
\end{eqnarray*}
Noting that $Q_{\theta(k)}\geq0, R_{\theta(k)}>0$, we have $u^{\ast}(k)=0, C_{\theta(k)}x^{\ast}(k)=0$.
That is,
\begin{eqnarray*}
x^{\ast}(k+1)=A_{\theta(k)}x^{\ast}(k), \ \ C_{\theta(k)}x^{\ast}(k)=0.
\end{eqnarray*}
Considering the exactly observable of $(C,A)$, it implies that $x=0$, which is a contradiction
with $x\neq0$. Hence, there must exist $N_{0}$, such that $P_{\theta(0)}^{N_{0}}(0)>0$.
Therefore,
\begin{eqnarray*}
P_{i}=\lim\limits_{N\rightarrow\infty}P_{i}^{N}(k)\geq P_{i}^{N_{0}}(k)>0.
\end{eqnarray*}

\begin{theorem}
Under Assumptions 1 and 2, the system (\ref{ff1}) is mean square stabilizable if and only if there exists a unique solution to (\ref{f40})-(\ref{f42}) such that $P_i>0, i=1,2,\cdots,L$. In this case, the controller
\begin{eqnarray}
u(k)=-\Upsilon_i^{-1}M_ix(k),k\geq 0\label{f43}
\end{eqnarray}
stabilizes (\ref{ff1}) in the mean square sense and minimizes the cost function (\ref{f38}). The optimal cost is given by
\begin{eqnarray}
J^{*}=\mbox{E}\{x_0'P_{\theta(0)}x_0\}.\label{f44}
\end{eqnarray}
\end{theorem}
\emph{Proof}.-\emph{Sufficiency}: Assume $P_i (i=1,2,\cdots,L)$ is a solution to (\ref{f40}) such that $P_i>0$. Firstly, we will show that (\ref{ff1}) is mean square stabilizable with the controller (\ref{f43}). For this purpose, define the Lyapunov function candidate $V(k,x(k))$ as
\begin{eqnarray}
V(k,x(k))=\mbox{E}[x(k)'P_{\theta(k)}x(k)].\label{f45}
\end{eqnarray}
The convergence of $V(k,x(k))$ is to be proven. Employing (\ref{ff1}) and (\ref{f40})-(\ref{f42}) yields
\begin{eqnarray}
&&V(k,x(k))-V(k+1,x(k+1))\nonumber\\
&=&\mbox{E}\{x(k)'P_{\theta(k)}x(k)-x(k+1)'P_{\theta(k+1)}x(k+1)\}\nonumber\\
&=&\mbox{E}\{x(k)'P_{\theta(k)}x(k)-[A_{\theta(k)}x(k)+B_{\theta(k)}u(k)]'P_{\theta(k+1)}[A_{\theta(k)}x(k)+B_{\theta(k)}u(k)]\}\nonumber\\
&=&\mbox{E}\{x(k)'Q_ix(k)+u(k)'R_iu(k)-[u(k)+\Upsilon_i^{-1}M_ix(k)]'\Upsilon_i[u(k)+\Upsilon_i^{-1}M_ix(k)]\}\label{f46}\\
&=&\mbox{E}\{x(k)'Q_{\theta(k)}x(k)+u(k)'R_{\theta(k)}u(k)\}\geq 0,k\geq 0, \label{f47}
\end{eqnarray}
where $u(k)=-\Upsilon_i^{-1}M_ix(k)$ for $k\geq 0$ has been used in (\ref{f46}). The above inequality (\ref{f47}) indicates that $V(k,x(k))$ decreases with respect to $k$. From Theorem 2, we know that
\begin{eqnarray}
V(k, x(k))=\mbox{E}\{x(k)'P_{\theta(k)}x(k)\}\geq 0,\label{f48}
\end{eqnarray}
which means that $V(k, (x(k)))$ is bounded, and thus is convergent.

Now let $m$ be any nonnegative integer. By adding from $k=m$ to $k=m+N$ on both sides of (\ref{f47}) and letting $m\rightarrow +\infty$, it yields that
\begin{eqnarray}
&&\lim_{m\rightarrow \infty}\sum_{k=m}^{m+N}\mbox{E}[x(k)'Q_{\theta(k)}x(k)+u(k)'R_{\theta(k)}u(k)]\nonumber\\
&=&\lim_{m\rightarrow \infty}[V(m, x(m))-V(m+N+1, x(m+N+1))]\nonumber\\
&=&0,\label{f49}
\end{eqnarray}
in which the last equality holds owning to the convergence of $V(k,x(k))$. Note that
\begin{eqnarray}
&&\sum_{k=0}^{N}\mbox{E}[x(k)'Q_{\theta(k)}x(k)+u(k)'R_{\theta(k)}u(k)]\geq \mbox{E}\{x_0'P_{\theta(0)}^N(0)x_0\}.\nonumber
\end{eqnarray}
Via a time-shift of length of $m$, it leads to
\begin{eqnarray}
&&\sum_{k=m}^{m+N}\mbox{E}[x(k)'Q_{\theta(k)}x(k)+u(k)'R_{\theta(k)}u(k)]\nonumber\\
&\geq& \mbox{E}\{x(m+0)'P_{\theta(m+0)}^{N+m}(m+0)x(m+0)\}\nonumber\\
&=&\mbox{E}\{x_m'P_{\theta(0)}^N(0)x_m\}\geq 0\label{f50}
\end{eqnarray}
In view of (\ref{f49}), we have proven that
\begin{eqnarray}
\lim_{m\rightarrow \infty}\mbox{E}\{x_m'P_{\theta(0)}^N(0)x_m\}=0, \forall N\geq 0.\label{f51}
\end{eqnarray}
 In the proof of Theorem 2, we have shown that there exists $N_0$, such that $P_{\theta(0)}^{N_{0}}(0)$ is positive definite. Thus (\ref{f51}) implies that $\lim_{m\rightarrow \infty}\mbox{E}[x_m'x_m]=0$. Therefore, the controller (\ref{f43}) stabilizes (\ref{ff1}) in the mean square sense.

Secondly, we will show that the cost function (\ref{f38}) is minimized by $(\ref{f43})$. Adding from $k=0$ to $k=N$ to (\ref{f46}) yields
\begin{eqnarray}
&&\mbox{E}\{\sum_{k=0}^N[x(k)'Q_{\theta(k)}x(k)+u(k)'R_{\theta(k)}u(k)]\}\nonumber\\
&=&V(0,x_0)-V(N+1,x(N+1))\nonumber\\
&&+\sum_{k=0}^N\mbox{E}\{[u(k)+\Upsilon_i^{-1}M_ix(k)]'\Upsilon_i[u(k)+\Upsilon_i^{-1}M_ix(k)]\},\label{f52}
\end{eqnarray}
in which $V(0,x_0)$ and $V(N+1,x(N+1))$ are defined in (\ref{f45}). Then $\lim_{k\rightarrow \infty}V(k,x(k))=0$ is to be shown.
%In fact, in view of (\ref{f50}), it follows that
%\begin{eqnarray}
%0\leq V(k,x(k))=\mbox{E}\{x(k)'P_{\theta(k)}x(k)\}.\nonumber
%\end{eqnarray}
Now we only consider the controller which stabilizes system (\ref{ff1}). Thus $\lim_{k\rightarrow \infty}\mbox{E}\{x(k)'P_{\theta(k)}x(k)\}=\lim_{k\rightarrow \infty}V(k,x(k))=0$. By letting $N\rightarrow \infty$ on both sides of (\ref{f52}), the cost function (\ref{f38}) is rewritten as \begin{eqnarray}
J=\mbox{E}\{x_0'P_{\theta(0)}x_0\}+\sum_{k=0}^{\infty}\mbox{E}\{[u(k)+\Upsilon_i^{-1}M_ix(k)]'\Upsilon_i[u(k)+\Upsilon_i^{-1}M_ix(k)]\},\label{f53}
\end{eqnarray}
In view of the positive definiteness of $\Upsilon_i, i=1,\cdots,L$, the optimal controller to minimizes (\ref{f53}) must be (\ref{f43}), and the corresponding optimal cost is as (\ref{f44}). Therefore, the proof of the sufficiency is finished.

\emph{Necessity}: Suppose the system (\ref{ff1}) is mean square stabilizable. In Theorem 2, the existence of the solution to (\ref{f40})-(\ref{f42}) satisfying $P_i>0 (i=1,2,\cdots,L)$ has been verified. We just need to show the uniqueness. Let $S_i (i=1,2,\cdots,L)$ be another solution to (\ref{f40})-(\ref{f42}) satisfying $S_i>0 (i=1,2,\cdots,L)$, i.e.,
\begin{eqnarray}
S_i&=&A_i'(\sum_{j=1}^L\lambda_{ij}S_j)A_i+Q_i-A_i'(\sum_{j=1}^L\lambda_{ij}S_j)B_i [B_i'(\sum_{j=1}^L\lambda_{ij}S_j)B_i+R_i]^{-1}\nonumber\\
&&\times B_i'(\sum_{j=1}^L\lambda_{ij}S_j)A_i,i=1,2,\cdots,L, \label{f54}\\
\Delta_i&=&B_i(\sum_{j=1}^L\lambda_{ij}S_j)B_i+R_i, i=1,2,\cdots,L, \label{f55}\\
\Pi_i&=&B_i(\sum_{j=1}^L\lambda_{ij}S_j)A_i, i=1,2,\cdots,L.\label{f56}
\end{eqnarray}
In view of the proof of sufficiency as in the above, the optimal value of the cost function (\ref{f38}) is as
\begin{eqnarray}
J^{*}=\mbox{E}\{x_0'P_{\theta(0)}x_0\}=\mbox{E}\{x_0'S_{\theta(0)}x_0\}.\nonumber
\end{eqnarray}
As $x_0$ is arbitrary,  the above equation implies that
$P_i=S_i, i=1,2,\cdots,L$. It follows from (\ref{f40})-(\ref{f42}), (\ref{f54})-(\ref{f56}) that, $\Upsilon_i=\Delta_i, M_i=\Pi_i (i=1,2,\cdots,L)$. Thus the uniqueness has been proven. The proof of necessity is now complete.

If $Q_i=I_{n\times n}, i=1,\cdots,L$ and $R_i=I_{m\times m}, i=1,\cdots,L$ in (\ref{f38}), where $I_{n\times n}$ is the identity matrix with dimension of $n$ and $I_{m\times m}$ is the identity matrix with dimension of $m$, the conditions of Assumption 1 and Assumption 2 are guaranteed naturally, and the performance index becomes as
\begin{eqnarray}
J&=&\mbox{E}\{\sum_{k=0}^{\infty}[x(k)'x(k)+u(k)'u(k)]\}.\label{ff38}
\end{eqnarray}
Then the stabilization solution to the infinite horizon problem (\ref{ff38}) can be stated as.
\begin{corollary}
The system (\ref{ff1}) is mean square stabilizable if and only if there exists a unique solution to (\ref{f40})-(\ref{f42}) such that $P_i>0, i=1,2,\cdots,L$.
In this case, the controller
\begin{eqnarray}
u(k)=-\Upsilon_i^{-1}M_ix(k),k\geq 0\label{f43}
\end{eqnarray}
stabilizes (\ref{ff1}) in the mean square sense and minimizes the cost function (\ref{f38}). The optimal cost is given by
\begin{eqnarray}
J^{*}=\mbox{E}\{x_0'P_{\theta(0)}x_0\}.\label{f44}
\end{eqnarray}
\end{corollary}

\begin{remark}
%(Comparison with the results in \cite{CV01},\cite{CV02})\\
In \cite{CV01}, a new detectability concept (weak detectability) for discrete-time MJLS was presented, and the new concept supplied a sufficient condition for the mean square stable for the infinite-horizon linear quadratic controlled system. The result can be summarized as:
If the system was weak detectable and there existed a positive semi definite solution to the CARE, then the system with the optimal feedback gain was mean square stable. Further, the necessary and sufficient conditions were supplied in \cite{CV02}, and we can summarize the result as: Under the assumption that the system was weak detectable, the system was mean square stabilizable if and only if there existed a positive semi-definite solution to the CARE. Although the sufficient conditions \cite{CV01} and necessary and sufficient conditions \cite{CV02} for the infinite-horizon stabilization problem were given. However, the computational test for weak detectability is not intuitive, and it is not easy to check. In Corollary 1, we give the necessary and sufficient conditions for the stabilization of the system without additional prerequisite. We just need to determine the existence of a positive definite solution to the CARE. It is easy to check.
\end{remark}
%%%%%%%%%%%%%%%%%%%%%%%%%%%%%%%%%%%%%%%%%%%%%%%%%%%%%%%%%%%%%%%%%%%%%%%%%%%%%%%%%%%%%%%%%%%%%%%%%%%%%%%%%%%%%%%%%%%%%%%%%%%%%%%%%%%%%%%%%%%%%%%
\section{Numerical Examples}
In this section, we present a simple example to illustrate the
previous theoretical results. Consider a second-order dynamic system (\ref{ff1}) with the performance ({\ref{f2}}). The specifications of the system and the weighting matrices are as follows
\begin{eqnarray}
&&A_1=\left[
      \begin{array}{cc}
        2 & 1.1 \\
        -1.7 & -0.8 \\
      \end{array}
    \right], A_2=\left[
                    \begin{array}{cc}
                      0.8 & 0 \\
                      0 & 0.6 \\
                    \end{array}
                  \right], B_1=\left[
                                  \begin{array}{c}
                                    1 \\
                                    1 \\
                                  \end{array}
                                \right],B_2=\left[
                                               \begin{array}{c}
                                                 2 \\
                                                 1 \\
                                               \end{array}
                                             \right],\nonumber\\
&&Q_{1}=\left[
          \begin{array}{cc}
            1& 0 \\
            0 & 1 \\
          \end{array}
        \right],Q_{2}=\left[
                          \begin{array}{cc}
                            1 & 0 \\
                            0 & 1 \\
                          \end{array}
                        \right],R_{1}=1,R_2=1.\nonumber
\end{eqnarray}
$\theta(k)$ is the Markov chain taking values in a finite set $\{1,\ 2\}$
with transition rate $\lambda_{11}=0.9$ and $\lambda_{22}=0.3$. The initial distribution of $\theta(k)$ is $(0.5,\ 0.5)$.
The initial
state $x(0)=[5\ 5]'$.

In this example, the time horizon is set to $N=20$. And the final penalty matrix $P_1(10)=I_2, P_2(10)=I_2$. Without loss of generality, we run $50$
Monte Carlo simulations from $k=0$ to $20$. The simulation results
are obtained as follows. Fig.\,1 shows a sample path of the Markov chain $\theta(k)\in \{1, 2\}$. The Riccati coefficients of the matrix $P_i(k)(i=1,2)$ obtained using MATLAB are shown in Fig.\,2. The optimal states are plotted in Fig.\,3
and the optimal control is shown in Fig.\,4.
% Finally, the sum of the three kinds of estimation error
%covariances of ${x}_1(k)$ and ${x}_2(k)$ are
%given in Fig.\,10 and Fig.\,11.
\begin{figure}[!h]
      \centering
      \includegraphics [scale=0.7]{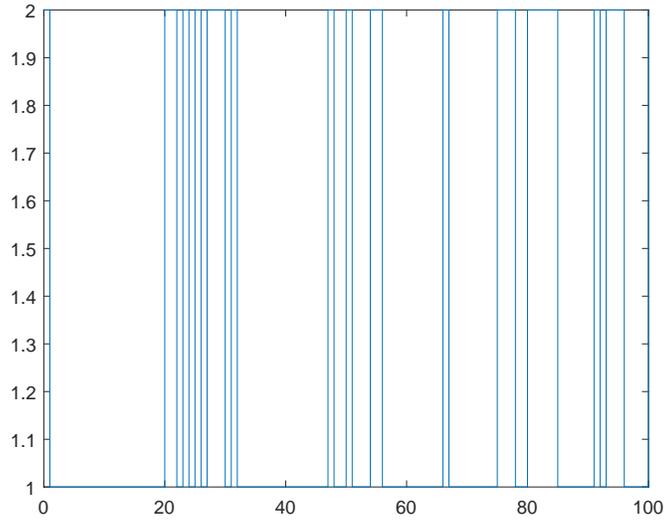}
      \caption{One path of the Markov chain $\theta(k)\in \{1,\ 2\}$}%over observations with time-varying delay}
      \label{f1}
\end{figure}
\begin{figure}[!h]
      \centering
      \includegraphics [scale=0.7]{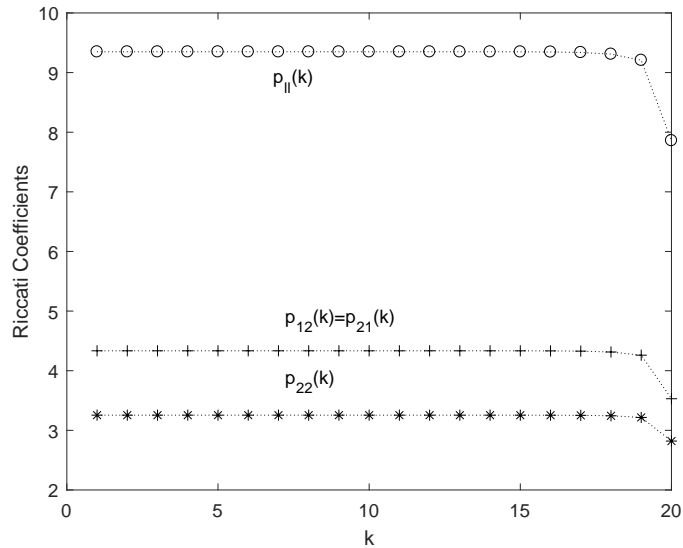}
      \caption{The Riccati coefficients of the matrix $P_i(k)(i=1,2)$}%over observations with time-varying delay}
      \label{f1}
\end{figure}
\begin{figure}[!h]
      \centering
      \includegraphics [scale=0.7]{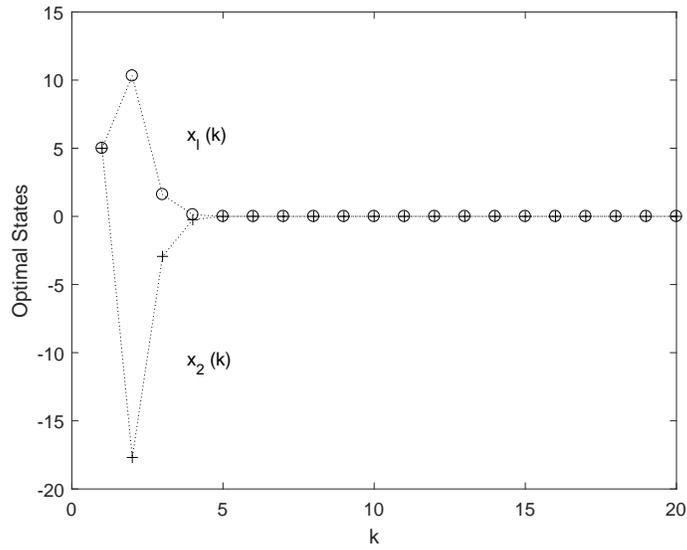}
      \caption{The optimal state trajectories}%over observations with time-varying delay}
      \label{f1}
\end{figure}
\begin{figure}[!h]
      \centering
      \includegraphics [scale=0.7]{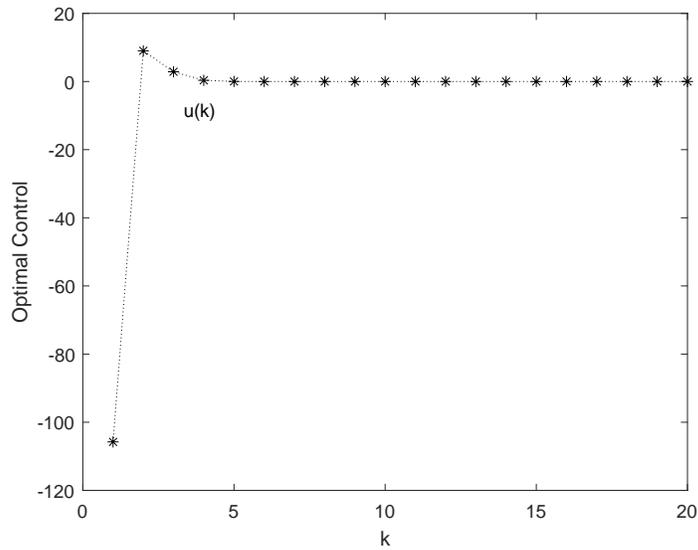}
      \caption{The optimal control}%over observations with time-varying delay}
      \label{f1}
\end{figure}
\section{Conclusions}
This paper has addressed the finite-horizon and infinite-horizon optimal control problems for the MJLS. A general situation in the former has been considered, and a necessary and sufficient condition for the existence of the optimal controller has been proposed for the first time. Later, we have proposed a necessary and sufficient condition for the mean square stabilizable of the MJLS. To show the existence of such a solution, one just need to prove the positiveness of the solution to the corresponding CARE. The condition is easily verifiable. As far as we know, no such conditions have been given for the mean mean square stabilizable of the MJLS before.% The result developed in this paper can be viewed as a generalization of the previous works developed in [18] for the multiplicative noise systems without delay.

\end{document}